\documentclass[12pt,reqno,a4wide]{amsart}

\usepackage{tikz} 
\usepackage{calrsfs}
\usepackage{mathrsfs}

\usepackage{array}
\def\P{\mathscr{P}}

\usepackage{hyperref}

\allowdisplaybreaks

\catcode`,\active

\catcode`\,12

\begin{document}
\bibliographystyle{plain}

%
%

	\title
	{Pell-Padovan tetranacci numbers and their Hadamard product with classical sequences}

	\author[H. Prodinger ]{Helmut Prodinger }
	\address{Department of Mathematics, University of Stellenbosch 7602, Stellenbosch, South Africa
	and
NITheCS (National Institute for
Theoretical and Computational Sciences), South Africa.}
	\email{warrenham33@gmail.com}


\begin{abstract}
A quick way to compute generating functions related to Pell-Padovan tetranacci numbers and classical sequences of recursions of order two is
provided. Eight special instances can be computed at once.
\end{abstract}

\maketitle

The sequence of \emph{Pell-Padovan tetranacci numbers}
\begin{equation*}
\P_{n+4}=\P_{n+2}+2\P_{n+1}+\P_{n},\quad \P_{0}=0,\;\P_{1}=1,\;\P_{2}=1,\;\P_{3}=1,
\end{equation*}
is the main interest in \cite{Zerr}.
The generating function is
\begin{equation*}
\sum_{n\ge0}\P_nz^n=\frac{z(1+z)}{1-z^2-2z^3-z^4}.
\end{equation*}
The Hadamard product of two generating functions is defined by
\begin{equation*}
\sum_{n\ge0}f_nz^n\odot\sum_{n\ge0}g_nz^n=\sum_{n\ge0}f_ng_nz^n.
\end{equation*}
The paper \cite{Zerr} computes
\begin{equation*}
	\sum_{n\ge0}\P_nz^n\odot\sum_{n\ge0}X_nz^n
\end{equation*}
for any of the 8 second order sequences in Table~\ref{la1}, by using a method with symmetric functions and discussing each instance separately.
We consider generally
\begin{equation*}
\frac{a+bz}{1+cd+dz^2}
\end{equation*}
for  $\sum_{n\ge0}X_nz^n$ and cover all 8 cases in one go. We obtain
\begin{equation*}
	\sum_{n\ge0}\P_nz^n\odot \frac{a+bz}{1+cd+dz^2} =\frac{\mathscr{N}}{\mathscr{D}}
\end{equation*}
with
\begin{multline*}
\mathscr{N}=	b-ac+  ( -ad-cb+a{c}^{2}  ) z+db{z}^{2}+d  ( ad-2 cb
	 ) {z}^{3}\\-d  ( -db+dac-{c}^{2}b  ) {z}^{4}+{d}^{2}
	 ( ad-cb  ) {z}^{5}
\end{multline*}
and
\begin{multline*}
\mathscr{D}=1+  ( -{c}^{2}+2 d  ) {z}^{2}+2 c  ( {c}^{2}-3 d
 ) {z}^{3}\\*+  ( -{c}^{4}+4 {c}^{2}d-{d}^{2}  ) {z}^{4}-
2 c{d}^{2}{z}^{5}+{d}^{2}  ( {c}^{2}+2 d  ) {z}^{6}-2 c{d}
^{3}{z}^{7}+{d}^{4}{z}^{8}.
\end{multline*}
The following Table~\ref{la1} shows how to specialize the parameters $a,b,c,d$ in order to obtain the 8 instances.

\begin{table}[h]
	\caption{Eight second order recursions}
	\label{la1}
\begin{tabular}{|l||l||c|c||c|c|c|c|}
	\hline
	  &recursion (term $=0$)&$X_0$&$X_1$ & $a$    & $b$&$c$&$d$ \\ \hline
	$k$-Fibonacci &$X_n-kX_{n-1}-X_{n-2}$& $1$&$k$&  1    & 0&$-k$&$-1$ \\
	$k$-Pell & $X_n-2X_{n-1}-kX_{n-2}$ & $0$&$1$&0    & 1&$-2$&$-k$ \\
		$k$-Jacobsthal&$X_n-kX_{n-1}-2X_{n-2}$ &   $0$&$1$& 0    & 1&$-k$&$-2$ \\
	$k$-Mersenne &$X_n-3kX_{n-1}+2X_{n-2}$&    $0$&$1$&0    & 1&$-3k$&$2$ \\
	Chebyshev first kind &$X_n-2xX_{n-1}+X_{n-2}$& $1$&$x$&   1    & $-x$&$-2x$&$1$ \\
	Chebyshev second kind &$X_n-2xX_{n-1}+X_{n-2}$& $1$&$2x$&   1    & $0$&$-2x$&$1$ \\
	Chebyshev third kind &$X_n-2xX_{n-1}+X_{n-2}$&   $1$&$2x-1$& 1    & $-1$&$-2x$&$1$ \\
	Chebyshev fourth kind &$X_n-2xX_{n-1}+X_{n-2}$&  $1$&$2x+1$&  1    & $1$&$-2x$&$1$ \\
	\hline
\end{tabular}

\end{table}
All these 8 sequences are covered in \cite{OEIS}.

Our method of choice is to use the Maple package \textsf{gfun}. It has various useful programs, like \textsf{`hadamardproduct`}, 
\textsf{`rec*rec`} etc. The easiest way to proceed seems to be to compute a list of, say, 30 coefficients of 
$\sum_{n\ge0}\P_nz^n\odot \frac{a+bz}{1+cd+dz^2}$ and use  \textsf{`guessgf`}; the program guesses than the generating function
$\mathscr{N}/\mathscr{D}$ as posted. This procedure is fully legitimate, as we know \emph{a priori} that  we will get a recursion of
order $4\times 2=8$. Knowing a sufficient number of initial coefficients, one knows that the guessed answer is indeed the correct answer. 
I learnt this type of argument from Doron Zeilberger.  Once one \emph{has} the rational generating function there are indeed various  ways
for  justification.

Of course, this approach is not restricted to the sequences discussed in this note.

\bibliographystyle{plain}

\end{document}